\documentclass[12]{amsart}

\usepackage{epsfig}
\usepackage{amsmath, amsthm, amssymb}
\usepackage{xcolor}
\usepackage{cancel}

\newtheorem{defn}{Definition}[section]

\newtheorem{thm}[defn]{Theorem}

\newtheorem{remark}[defn]{Remark}

\newcommand{\be}{\begin{equation}}
\newcommand{\ee}{\end{equation}}
\newcommand{\bea}{\begin{eqnarray}}
\newcommand{\eea}{\end{eqnarray}}
\newcommand{\beas}{\begin{eqnarray*}}
\newcommand{\eeas}{\end{eqnarray*}}

\newcommand{\R}{\mathbb{R}}

\newcommand\Prb{\mathbb{P}}
\newcommand\E{\mathbb{E}}

\newcommand{\goto}{\rightarrow}
\newcommand{\hsp}{\hspace{.3in}}
\newcommand{\bp}{\begin{proof}}
\newcommand{\ep}{\end{proof}}

\begin{document}

\title[Gelation and Extinction]{Gelation in Vector Multiplicative Coalescence and Extinction in Multi-Type Poisson Branching Processes}

\author{Heshan Aravinda}
\address{Department of Mathematics \& Statistics, University of Wyoming, Laramie, WY  82071, USA}
\email{hpathira@uwyo.edu}

\author{Yevgeniy Kovchegov}
\address{Department of Mathematics, Oregon State University, Corvallis, OR  97331, USA}
\email{kovchegy@math.oregonstate.edu}

\author{Peter T. Otto}
\address{Department of Mathematics, Willamette University, Salem, OR 97302, USA}
\email{potto@willamette.edu}

\author{Amites Sarkar}
\address{Department of Mathematics, Western Washington  University, Bellingham, WA  98225, USA}
\email{amites.sarkar@wwu.edu}


\date{}



\begin{abstract}

In this note, we present a novel connection between a multi-type (vector) multiplicative coalescent process and a multi-type branching process with Poisson offspring distributions.  More specifically, we show that the equations that govern the phenomenon of gelation in the vector multiplicative coalescent process are equivalent to the equations that yield the extinction probabilities of the corresponding multi-type Poisson branching process.  We then leverage this connection with two applications, one in each direction.  The first is a new quick proof of gelation in the vector multiplicative coalescent process, and the second is a new series expression for the extinction probabilities of the multi-type Poisson branching process. We also use random graphs to give a new derivation of the solution to the modified Smoluchowski coagulation equations, which describe the vector multiplicative coalescent process.
\end{abstract}

\maketitle
\section{Introduction}
The main motivation for this work is paper \cite{KO} where the authors used the multidimensional Lambert-Euler inversion in the analysis of the vector multiplicative coalescent process. More specifically, they solved the multidimensional Lambert-Euler inversion and utilized the solution to investigate the gelation phenomenon in the coagulation model with the vector-multiplicative kernel. Our goal in this note is to present an interesting connection between the vector multiplicative coalescent process and a multi-type branching process with Poisson offspring distributions. As a consequence, we give a quick proof of the gelation phenomenon. Additionally, this connection leads to a series expression for the extinction probabilities of the multi-type Poisson branching process. This illuminates a connection between these two paradigms which we hope will continue to reveal new insights and potential for cross-disciplinary research.\vskip2mm
We begin by recalling some of the notation and concepts from the theory of vector multiplicative coalescent processes. For a detailed account of this theory, we refer the reader to \cite{KOY, KO}.
\subsection{Vector multiplicative coalescent processes}
\label{vec_coal}
Roughly speaking, a vector multiplicative coalescent process is a Markov process that describes cluster merger dynamics, where the weight of each cluster is a finite-dimensional vector. Following the convention in \cite{KO}, to define this process formally, we use bra-ket notation. Specifically, $|{\bf x}\rangle$ will denote the column vector representation of vector ${\bf x} \in \mathbb{R}^k$,
and $\langle {\bf x}|$ will denote the row vector representation of vector ${\bf x} \in \mathbb{R}^k$.
For $c \in \mathbb{R}$ and ${\bf x} \in \mathbb{R}^k$, $\,c|{\bf x}\rangle$ will represent the product $c{\bf x}$, a column vector.
Correspondingly, $\langle {\bf x}|{\bf y}\rangle=\langle {\bf y}|{\bf x}\rangle$ will be the dot product of ${\bf x}$ and ${\bf y}$ in $\mathbb{R}^k$.
Finally, for a matrix $M \in \mathbb{R}^{k \times k}$, $\langle {\bf x}|M|{\bf y}\rangle$ will represent the product ${\bf x}^{\sf T}M{\bf y}$ resulting in a scalar. \medskip

Consider a system with $k$ types of particles: $1,\hdots,k$.
For a given vector ${\boldsymbol \alpha} = (\alpha_1, \alpha_2, \ldots, \alpha_k) \in(0,\infty)^k$, the process begins
with $\langle {\boldsymbol \alpha} n | {\boldsymbol 1} \rangle$ singletons distributed between the $k$ types so that for each $i$, there are $\lfloor\alpha_i\,n\rfloor$ particles of type $i$.
From now on, we omit floor functions for the sake of clarity.

\medskip

Let $V \in \mathbb{R}^{k \times k}$ be a nonnegative, irreducible and symmetric matrix, which we will refer to as the {\it partition interaction matrix}. A particle of type $i$ bonds with a particle of type $j$ with the intensity rate $v_{i,j}/n$, where $v_{i,j}$ is the $(i,j)$ element in the matrix $V$. The bonds are formed independently. This process is called the {\it vector-multiplicative coalescent process}. Note that the irreducibility condition on $V$ ensures that particles of each type could ultimately coalesce with particles of every other type at some time in the vector-multiplicative coalescent process.\medskip

The {\it size} of a cluster of particles is represented by a $k$-dimensional vector ${\bf x} = (x_1, x_2, \ldots, x_k) \in \mathbb{Z}_{\geq 0}^k$ such that $\sum_{i=1}^k x_i>0$, where each $x_i$ represents the number of particles of type $i$
in the cluster. In other words, each cluster of weight ${\bf x}$ bonds together $x_1,x_2,\ldots, x_k$  particles of corresponding types $1,2,\ldots,k$. At any given time $t$, each pair of clusters with respective size vectors ${\bf x}$ and ${\bf y}$ will coalesce into a cluster of size ${\bf x+y}$ with rate $K({\bf x},{\bf y})/n$,
where
\be\label{VMkernel}
K({\bf x}, {\bf y}) = \langle {\bf x} | V | {\bf y} \rangle.
\ee
The last merger will create a cluster of size $\,{\boldsymbol \alpha} n$. Clusters coalesce independently of other clusters. Note that \eqref{VMkernel} can be interpreted at the level of individual particles: two clusters bond whenever particles from that cluster bond. This interpretation enables one to make the connection between coalescent processes and random graphs, which will be important later in the paper.\medskip

The kernel $K({\bf x}, {\bf y})$ defined in \eqref{VMkernel} will be referred to as the {\it vector-multiplicative kernel}.
The kernel is symmetric
$$K({\bf x}, {\bf y}) =K({\bf y},{\bf x}) \quad \text{ for all vectors }~{\bf x}, {\bf y}$$
and bilinear
$$K(c_1{\bf x} + c_2{\bf y}, {\bf z}) = c_1\,K({\bf x}, {\bf z}) + c_2\,K({\bf y}, {\bf z}) \quad \text{ for all vectors }~{\bf x}, {\bf y}, {\bf z} ~~\text{ and scalars }~c_1,c_2.$$

In the vector-multiplicative coalescent process, let $\zeta_{\bf x}^{[n]}(t)$ denote the number of clusters of size ${\bf x}$ at time $t \geq 0$. The initial values are $\zeta_{\bf x}^{[n]}(0)=\alpha_i n$ if ${\bf x} = {\bf e}_i$, where ${\bf e}_i$ is the $i$-th standard basis vector, and $\zeta_{\bf x}^{[n]}(0)=0$ otherwise.
The process
$${\bf ML}_n(t)=\Big(\zeta_{\bf x}^{[n]}(t) \Big)_{{\bf x} \in \mathbb{Z}_{\geq 0}^k : \langle {\bf x}|{\bf 1}\rangle>0}$$
that counts clusters of all sizes in the vector-multiplicative coalescent process is the corresponding
{\it Marcus-Lushnikov process}.
For fixed ${\bf x}$ and $t$, we have that $\lim_{n \goto \infty} \frac{\zeta_{\bf x}^{[n]}(t)}{n}=\zeta_{\bf x}(t) $, where the convergence is in probability.
Thus $\zeta_{\bf x}(t)$ is the deterministic limiting fraction of clusters of each size ${\bf x}$. The differential equations that govern the dynamics of $\zeta_{\bf x}(t)$ are called the {\it Smoluchowski coagulation system of equations} and are given by
\be
\label{eqn:SmolEqn}
{d \over dt}\zeta_{\bf x}(t) =  - \zeta_{\bf x} \sum_y \zeta_{\bf y} \langle {\bf x} | V | {\bf y} \rangle  + \frac{1}{2}\sum_{{\bf y}, {\bf z} \,: {\bf y} + {\bf z} = {\bf x}} \langle {\bf y} | V | {\bf z} \rangle \zeta_{{\bf y}} \zeta_{{\bf z}}
\ee
with initial conditions $\zeta_{\bf x}(0)= \alpha_i$ if $x = {\bf e}_i$, and $\zeta_{\bf x}(0)= 0$ otherwise.  These equations describe the evolution of the size distribution of clusters in a system where clusters merge multiplicatively as defined in \eqref{VMkernel}. More specifically, the first term on the right hand side in \eqref{eqn:SmolEqn} corresponds to the depletion of clusters of size {\bf x} after merging with other clusters, and the second term corresponds to coagulation of smaller clusters to produce one of size {\bf x}. For further details on the Smoluchowski coagulation equations in various settings, including applications, we refer the reader to \cite{S,No,HP,B,K}. \bigskip

\vfill\eject

\textbf{Gelation.} For many kernels, including the multiplicative kernel of this paper, the coalescent process undergoes a phase transition known as \textit{gelation}. Informally, at a certain finite and deterministic time $T_{gel}$, a \textit{gel} forms, containing, as $n\to\infty$, a positive proportion of all the particles. Gelation can only be observed ``as $n\to\infty$", since, for fixed $n$, the number of particles is finite, and so every cluster contains a positive fraction of the total number of particles. Mathematically, at $T_{gel}$, two things happen simultaneously~\cite{KO}:

\medskip

$\bullet$ The second moments of $\zeta_{\bf x}(t)$, i.e., the entries in the matrix $\sum_{\bf x}\zeta_{\bf x}(t)|{\bf x}\rangle\langle {\bf x}|$, are no longer finite.

\smallskip

$\bullet$ The limiting total mass $\,\sum_{\bf x} \zeta_{\bf x}(t) |{\bf x} \rangle \,$ is no longer conserved.

\medskip

Of course, for fixed $n$, mass \textit{is} conserved: gelation can only be observed in the limit $n\to\infty$. But, after $T_{gel}$, the \textit{limiting} mass is not conserved because the gel, invisible in \eqref{eqn:SmolEqn}, has swallowed some positive fraction of finite-size clusters, and only these clusters are modelled in \eqref{eqn:SmolEqn}. From now on, we will simply say that mass is not conserved after gelation, and omit reference to the underlying limiting process. For more details, see~\cite{J,R,AIM,VE,HEZ}. Gelation corresponds exactly to the formation of a giant component in an associated random graph, which we describe later.

\subsection{The solution to the Smoluchowski equations prior to gelation} In \cite{KO}, it was shown that, prior to gelation, the solutions to the Smoluchowski equations and the following {\it modified Smoluchowski equations (MSE)} given by
\be
\label{eqn:VMSmolEqn}
{d \over dt}\zeta_{\bf x}(t) =  - \zeta_{\bf x} \langle {\bf x} | V | \boldsymbol{\alpha} \rangle  + \frac{1}{2}\sum_{{\bf y}, {\bf z} \,: {\bf y} + {\bf z} = {\bf x}} \langle {\bf y} | V | {\bf z} \rangle \zeta_{{\bf y}} \zeta_{{\bf z}}
\ee
with the same initial conditions, coincide. This is because, prior to gelation, mass is conserved, and so
\[
\sum_{\bf x} \zeta_{\bf x}(t) |{\bf x} \rangle  =  |\boldsymbol{\alpha} \rangle.
\]
Combining this with \eqref{eqn:SmolEqn} yields \eqref{eqn:VMSmolEqn}.

The unique solution $\zeta_{\bf x}(t)$ of \eqref{eqn:VMSmolEqn} was derived in \cite{KO} (see Corollary 3.11).
Specifically, for a vector ${\bf x} \in \mathbb{Z}_{\geq 0}^k$, let
$\,{\bf x}!=x_1!x_2!\hdots x_k!\,$,
and for vectors ${\bf a}$ and ${\bf b}$ in $\mathbb{R}_{+}^k$, let
$\,{\bf a}^{\bf b}=a_1^{b_1} a_2^{b_2}\hdots a_k^{b_k}$.
Now, consider a complete graph $K_k$ consisting of vertices $\{1,\hdots,k\}$ with weights $w_{i,j}=w_{j,i}\geq 0$ assigned to
its edges $[i,j]$ ($i \not=j$).
Let the weight $W(\mathcal{T})$ of a spanning tree $\mathcal{T}$ be the product of the weights of all of its edges.
Finally, let $\tau(K_k, w_{i,j})=\sum\limits_{\mathcal{T}}W(\mathcal{T})$ denote the {\it weighted spanning tree enumerator},
i.e., the sum of weights of all spanning trees in $K_k$.  In the case $w_{i,j} \equiv 1$, the weighted spanning tree enumerator simply equals the number of spanning trees of the graph.
In terms of these notations, the solution to the modified Smoluchowski equations that yield the limiting fraction of clusters of size ${\bf x}$ is given by
\be\label{eqn:MSEsol}
\zeta_{\bf x}(t) = {1 \over {\bf x}!}\boldsymbol{\alpha}^{\bf x} {\tau(K_k, x_i x_j v_{i,j}) \over {\bf x}^{\bf 1}}(V {\bf x})^{{\bf x}- {\bf 1}} e^{-\langle {\bf x} | V | \boldsymbol{\alpha} \rangle t} t^{\langle {\bf x} |\boldsymbol{1} \rangle -1}.
\ee
\vskip2mm \noindent

\textbf{Outline of the paper.\,} The rest of the paper is organized as follows: Section 2 is devoted to a short proof of equation \eqref{eqn:MSEsol}, based on the theory of random graphs. It provides a random graphs perspective on the Smoluchowski equations, which we hope is of independent interest. In Section 3, we discuss the connection between gelation and multidimensional Lambert-Euler inversion. In Section 4, we recall some basic concepts from the theory of multi-type branching processes. We conclude in Section 5 by stating the connection between the vector-multiplicative coalescent process and multi-type Poisson branching processes, and giving two applications of this connection.

\section{Random graphs and an alternative proof of equation \eqref{eqn:MSEsol}}
The following proof uses the connection between the vector multiplicative coalescent process and the theory of {\it random graphs}.
Random graphs model networks with random connections. They were introduced by Erd\H os and R\'enyi in 1960~\cite{ER}. The basic model $G(n,p)$ has $n$ vertices, with edges inserted independently at random with probability $p$. One way to think of this model is to imagine $n$ fixed and very large, with $p$ slowly increasing. Specifically, we will set $p=t/n$, where $t$ is a time parameter, so that, as time progresses, the graph $G=G(n,t/n)$ gathers more and more edges, progressing from a collection of isolated vertices to the complete graph with all edges present. Along the way, $G$ undergoes several {\it phase transitions}. One of the most interesting
is the emergence of the so-called {\it giant component} at $p=1/n$, i.e., at time $t=1$. Roughly speaking, at fixed time $t=1-\epsilon$, as $n\to\infty$, all the connected components of $G$ have size $O(\log n)$, but
at time $t=1+\epsilon$, the largest connected component of $G$ (the ``giant") has size $C(\epsilon)n$. As with gelation, the emergence of the giant component is an asymptotic phenomenon, which can only be observed in the limit as $n\to\infty$. It was studied in detail by Erd\H os and R\'enyi~\cite{ER}; subsequent work by Bollob\'as and others led to the field of random graphs~\cite{Bol01}.\medskip

The connection between the vector multiplicative coalescent process $C$ and (an inhomogeneous version of) the above random graph $G$ works as follows. Small clusters of particles in $C$ correspond to small connected components of $G$. The giant component of $G$ corresponds to the gel in $C$. The correspondence between $C$ and $G$ was first described in \cite{A} for the Erd\H os-R\'enyi model (which corresponds to a coalescent process with only one type of particle), and in~\cite[Section 1.2]{KO} for the vector multiplicative process with $k$ particle types. We refer to~\cite{KO} for a full account. Here, we include only enough detail to describe the proof of \eqref{eqn:MSEsol}.\medskip

Our random graph $G$, whose $n$ vertices are partitioned into $k$ sets $A_1,\ldots, A_k$, models the coalescent process at a fixed time $t$. The $n$ particles in $C$ correspond to the $n$ vertices of $G$, and, as before, we will let $n\to\infty$. The numbers $\alpha_i n$ of particles of each type in $C$ correspond to the sizes of the parts $A_i$ of $G$, so that $G$ consists of $n$ vertices, divided into $k$ parts $A_1,\ldots,A_k$, in such a way that part $A_i$ contains $\alpha_i n$ vertices. Each potential edge $x_ix_j$ of $G$ (with $x_i\in A_i$ and $x_j\in A_j$)
is included in $E(G)$, the edge set of $G$, with probability
\[
1-e^{-v_{ij}t/n}=\frac{v_{ij}t}{n}+O\left(\left(\frac{v_{ij}t}{n}\right)^2\right).
\]
Here, $V=(v_{ij})$ is the partition interaction strength matrix defined above. The exponential formulation is convenient (but not strictly necessary) so that $G$ can be seen as a snapshot of a certain
memoryless process. Finally, for our purposes, $t$ will be a constant, not depending on $n$. Note that if $k=\alpha_1=v_{11}=1$ then $G$ is (modulo the approximation above) just the Erd\H os-R\'enyi graph $G(n,t/n)$, with average degree (approximately) $t$. (Also note that we are allowing edges inside the parts, so that our graphs are not necessarily $k$-partite.)\medskip

With ${\boldsymbol\alpha}$ and $V$ fixed, there will be a critical value $t'$ of $t$ below which there is almost surely no giant component, and above which there is almost surely a giant component. Denote by $D[{\boldsymbol\alpha}]$ the diagonal matrix with the $\alpha_i$ on the diagonal. Then, this critical value $t'$ is such that the largest eigenvalue of the matrix $t'V D[{\boldsymbol\alpha}]$ equals 1 (i.e., $t'$ is the reciprocal of the spectral radius of $VD[{\boldsymbol\alpha}]$)~\cite{KO}.\medskip

Next, again with $V$ and ${\boldsymbol\alpha}$ fixed, and for given values of ${\bf x}=(x_1,\ldots,x_k)$ and $t$, we seek the scaled number of {\it finite} components in $G$ containing $x_i$ vertices of type $i$. By ``scaled", we mean the number of such components divided by $n$. Call (the limit as $n\to\infty$ of) this scaled quantity $\phi_{\bf x}(t)$. Note that the process undergoes a qualitative change when $t=t'$, but, for a given ${\bf x}$, the functions $\phi_{\bf x}(t)$ are in fact continuous for all $t$~\cite{Bol01,KO}. Informally, when it comes to counting {\it components} instead of vertices, the giant component is asymptotically irrelevant, since there is only one of it, and we are dividing by $n$.\medskip

The following theorem, connecting the process $C$ and the graph $G$, appears in \cite{KO} (see also \cite{A}).

\begin{thm}\cite{KO}
For given ${\boldsymbol\alpha},V$ and for any ${\bf x}$ and $t<t'=T_{\rm gel}$, we have
\[
\phi_{\bf x}(t)=\zeta_{\bf x}(t),
\]
where the functions $\zeta_{\bf x}(t)$ are the solutions of the modified Smoluchowski equations \eqref{eqn:VMSmolEqn} above.
\end{thm}

Consequently, at least in the pre-gelation interval, we can use the random graph formulation to compute $\zeta_{\bf x}(t)$. We choose a method described to prove Theorem 11.4.2 in Chapter 11 of~\cite{AS}.\medskip

For simplicity, we first address the
case where $k=\alpha_1=v_{11}=1$. In this case, setting ${\bf x}=l$, we have
\be \label{eqn:RGMSEsol}
\zeta_l(t)=\phi_l(t)=\frac{l^{l-2}t^{l-1}e^{-lt}}{l!}.
\ee

The random graphs proof of this is as follows. We calculate $p_l(t)=l\phi_l(t)$, which is the probability that a given vertex
lies in a connected component $G_l$ of size $l$ (which, with probability tending to 1, is a tree - see~\cite{Bol01} for a proof). Now, setting $p=t/n$ for the
edge probability, we see that:
\begin{itemize}
\item There are $\binom{n-1}{l-1}$ choices for the other vertices in $G_l$.
\item There are $l^{l-2}$ labeled trees on these chosen vertices.
\item Each tree is present with probability $p^{l-1}(1-p)^{\binom{l}{2}-l+1}$.
\item The edges between $G_l$ and the other vertices are missing with probability $(1-p)^{l(n-l)}$.
\end{itemize}
Therefore, as $n\to\infty$,
\begin{align*}
p_l(t)&\sim\binom{n-1}{l-1}\cdot l^{l-2}\cdot p^{l-1}\cdot (1-p)^{l(n-l)}
\sim\frac{n^{l-1}}{(l-1)!}\cdot l^{l-2}\cdot \left(\frac{t}{n}\right)^{l-1}\cdot e^{-pln}\\
&=\frac{1}{(l-1)!}\cdot l^{l-2}\cdot t^{l-1}\cdot e^{-tl}
=\frac{l^{l-1}t^{l-1}e^{-lt}}{l!},
\end{align*}
which, recalling that $p_l(t)=l\phi_l(t)$, establishes \eqref{eqn:RGMSEsol}.

\medskip

Next we consider the general case. For this we once again require the {\it weighted spanning tree enumerator} defined above,
although here it arises in a slightly different form. Write $K_{\bf x}(V)$ for the complete graph on $\langle {\bf x} |\boldsymbol{1} \rangle$ vertices with
$x_i$ vertices of type $i$, and where the weight of an edge between vertices of types $i$ and $j$ is $v_{ij}$. Suppressing the dependence
on $V$, write $T_{\bf x}$ for the sum of the weights of all the spanning trees of $K_{\bf x}(V)$. Then, from Theorem 3.8 in~\cite{KO},
\[
T_{\bf x}=\frac{\tau(K_k, x_i x_j v_{i,j})}{ {\bf x}^{\bf 1}} (V {\bf x})^{{\bf x}- {\bf 1}}.
\]

Now, generalizing the argument above, pick a random vertex $v$.
The probability that $v$ is of type $i$ is $\alpha_i$. For simplicity of notation, suppose that $i=1$. We aim to compute
$p_{\bf x}(t)$, which is the probability that $v$ lies in a connected component $G_{\bf x}$ of type ${\bf x}$. Again, with probability
tending to 1, $G_{\bf x}$ is a tree~\cite{Bol01}. Omitting floor and ceiling functions:
\begin{itemize}
\item There are $\binom{\alpha_1n}{x_1-1}\binom{\alpha_2n}{x_2}\cdots\binom{\alpha_kn}{x_k}$ choices for the other vertices in $G_{\bf x}$
\item These form a tree with probability asymptotically $T_{\bf x}\cdot(t/n)^{\langle {\bf x} |\boldsymbol{1} \rangle -1}$
\item The edges from $G_{\bf x}$ to $G\setminus G_{\bf x}$ are missing with probability $\sim\exp(-\sum_i\sum_jx_iv_{ij}\alpha_j t)$
\end{itemize}
Therefore, as $n\to\infty$,
\begin{align*}
p_{\bf x}(t)&\sim \langle {\bf x} |\boldsymbol{1} \rangle \frac{n^{\langle {\bf x} |\boldsymbol{1} \rangle -1}\alpha_1^{x_1}\cdots\alpha_k^{x_k}}{x_1!\cdots x_k!}
\cdot T_{\bf x}
\cdot\left(\frac{t}{n}\right)^{\langle {\bf x} |\boldsymbol{1} \rangle -1}
\cdot{\rm exp}\left(-\sum_i\sum_j x_iv_{ij}\alpha_jt\right)\\
&= \langle {\bf x} |\boldsymbol{1} \rangle \frac{\alpha_1^{x_1}\cdots\alpha_k^{x_k}}{x_1!\cdots x_k!}
\cdot T_{\bf x}
\cdot t^{\langle {\bf x} |\boldsymbol{1} \rangle -1}
\cdot{\rm exp}\left(-\sum_i\sum_j x_iv_{ij}\alpha_jt\right),
\end{align*}
so that, since
\[
p_{\bf x}(t)=\phi_{\bf x}(t) \langle {\bf x} |\boldsymbol{1} \rangle,
\]
we have
\[
\zeta_{\bf x}(t)=\phi_{\bf x}(t)=\frac{\alpha_1^{x_1}\cdots\alpha_k^{x_k}}{x_1!\cdots x_k!}
\cdot T_{\bf x}
\cdot t^{\langle {\bf x} |\boldsymbol{1} \rangle -1}
\cdot{\rm exp}\left(-\sum_i\sum_j x_iv_{ij}\alpha_jt\right),
\]
proving \eqref{eqn:MSEsol}.

\section{Gelation and Lamber-Euler Inversion} \label{GelationLE}

For the coalescent process defined in Section \ref{vec_coal}, the initial total mass vector is assumed to be $\,\sum_{\bf x} \zeta_{\bf x}(0) |{\bf x} \rangle= |\boldsymbol{\alpha} \rangle$. As already discussed, we define gelation as the loss of total mass $\,\sum_{\bf x} \zeta_{\bf x}(t) |{\bf x} \rangle \,$ of the system after a critical time $T_{gel}$ called the \textit{gelation time}. That is
 \be
 \label{eqn:GelTimeDefn}
 T_{gel}=\inf\Big\{t>0\,:\,\sum_{\bf x} \zeta_{\bf x}(t) |{\bf x} \rangle  <  |\boldsymbol{\alpha} \rangle \Big\}.
 \ee
 In \cite[Corollary 4.1]{KO}, the existence of gelation was proved for the vector multiplicative coalescent process and the gelation time was given by
 \be
 \label{eqn:GelTime}
 T_{gel}={1 \over \rho ( V D[{\boldsymbol\alpha}] )},
 \ee
where $\rho(A)$ denotes the spectral radius (i.e., the largest of the absolute values of the eigenvalues) of the matrix $A$.\medskip

 The existence of gelation for the vector multiplicative coalescent process was proved by solving a multidimensional Lambert-Euler inversion problem.  This is a higher dimensional generalization of the equation originally studied by Lambert and followed up by Euler that gave rise to the well known Lambert W function \cite{JL, EL}. For a given vector $\boldsymbol{\alpha}\in(0,\infty)^k$, consider the region
\begin{equation}\label{eqn:R0}
R_0 =\left\{\boldsymbol{\alpha}  \in (0,\infty)^k :\, \rho\left(V D[{\boldsymbol\alpha}]\right) <1 \right\},
\end{equation}
its closure within $(0,\infty)^k$,
\begin{equation}\label{eqn:R0bar}
\overline{R}_0=\left\{\boldsymbol{\alpha}  \in (0,\infty)^k :\, \rho\left(V D[{\boldsymbol\alpha}] \right) \leq 1\right\},
\end{equation}
and the complement of $\overline{R}_0$ within $(0,\infty)^k$,
\begin{equation}\label{eqn:R1}
R_1 =\left\{\boldsymbol{\alpha} \in (0,\infty)^k :\, \rho\left(V D[{\boldsymbol\alpha}] \right) >1 \right\}.
\end{equation}

The following theorem implies the existence of a unique solution to the multidimensional Lambert-Euler problem (see \cite[Theorem 1.1]{KO}).
\begin{thm}[Multidimensional Lambert-Euler inversion]\label{thm:multEulerLambert}
Consider a nonnegative irreducible symmetric matrix $V \in \mathbb{R}^{k \times k}$. Then, for any given $\boldsymbol{\alpha} \in(0,\infty)^k$ and $t \geq 0$, there exists a unique vector ${\bf y}\in \overline{R}_0$ such that
\begin{equation}\label{eqn:yzSolsCoord}
y_i  e^{-\langle {\bf e}_i  | V | {\bf y} \rangle} = \alpha_i t  e^{-\langle {\bf e}_i  | V | \boldsymbol{\alpha} t \rangle} \qquad i=1,\hdots,k.
\end{equation}
Moreover, if $\boldsymbol{\alpha} t \in \overline{R}_0$, then ${\bf y}= \boldsymbol{\alpha} t$. If $\boldsymbol{\alpha} t \in R_1$,
then ${\bf y}< \boldsymbol{\alpha} t$ ($y_i < \alpha_i t$ for all $i$), i.e., ${\bf y}$ is the smallest solution of \eqref{eqn:yzSolsCoord}.
\end{thm}
Finally, in \cite[Lemma 4.2]{KO}, it was shown that the total mass $\sum_{\bf x}  \zeta_{\bf x}(t)|{\bf x}\rangle$ of the vector multiplicative coalescent process could be expressed in terms of the solution of the multidimensional Lambert-Euler inversion problem as follows:
\be
\label{eqn:TotalMass}
\sum_{\bf x}  \zeta_{\bf x}(t)|{\bf x}\rangle = {1 \over t}|{\bf y}(t)\rangle.
\ee
Then, Theorem \ref{thm:multEulerLambert} proves the existence of the gelation phenomenon for the vector multiplicative coalescent process and the gelation time given in \eqref{eqn:GelTime}.  Specifically,
\[
\sum_{\bf x}  \zeta_{\bf x}(t)|{\bf x}\rangle  \left\{ \begin{array}{ll}
= |{\boldsymbol \alpha}\rangle & \mbox{if $t \leq 1/\rho\left(V D[{\boldsymbol\alpha}] \right)$} \\
< |{\boldsymbol \alpha}\rangle & \mbox{if $t > 1/\rho\left(V D[{\boldsymbol\alpha}] \right)$}.
\end{array}
\right.
\]

\noindent{\bf Remark.} The one dimensional Lambert-Euler inversion reads as
    \begin{align*}
    \label{eqn: LE}
x(t) = \min\{x\geq 0\,:\, xe^{-x}= te^{-t}\}.
  \end{align*}
    In their now famous paper \cite{ER}, Erd\H{o}s and R\'enyi used the function $x(t)$ to establish the formation of a giant component in the theory of random graphs. In fact, in the random graph setting, the quantity $x(t)$ can be interpreted as \textit{the average degree outside the giant component} in an Erd\H{o}s-R\'enyi graph $G(n,t/n)$. It would be interesting to see if (\ref{eqn:yzSolsCoord}) can be used to obtain refined information about the vertices in the giant component of a random inhomogeneous graph, in the same spirit as the classical case. \vskip2mm
   Note that in this one-dimensional setting, $x(t)=t$ on $(0,1]$, and $x(t)<t$ on $(1, \infty)$. Moreover, $x(t) \downarrow 0$ monotonically as $t \to \infty$. On the other hand, the post-gelation behavior of the multi-dimensional solution to (\ref{eqn:yzSolsCoord}) can differ from the one-dimensional case. For example, consider the bipartite case with $V = \begin{vmatrix}
       0 & 1\\ 1 & 0
   \end{vmatrix}$ and ${\boldsymbol{\alpha}}= (15,2)$. Then, the critical time is given by $T_{gel} = \frac{1}{\sqrt{30}}$, and for $t> T_{gel}$, the graph of $y(t) = (y_1(t), y_2(t))$ is as follows.
   \begin{figure}[h]
    \centering
    \includegraphics[width=0.25\textwidth]{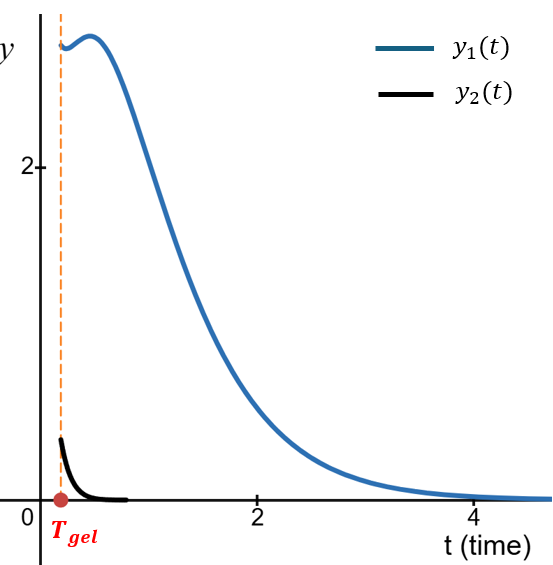}
    \caption{Post-gelation behavior of Lambert-Euler inversion}
    \label{fig:postge;}
\end{figure}
\vskip2mm
Observe that the function $y_1(t)$ decreases, increases, and finally decreases as $t$ increases. Numerical simulation seems to imply that there is a critical value for the ratio of the parameters in ${\boldsymbol{\alpha}}$ after which the solutions $(y_1(t), y_2(t))$ are no longer monotonically decreasing. Note, however, that only $(y_1(t)+y_2(t))/t$ has physical significance (as the total mass outside the gel).

\section{Multi-Type Poisson Branching Processes}

In this and the next sections, we give an overview of the relevant theory of multi-type branching processes, which is a generalization of the classical Galton-Watson process, extending it to multiple types of individuals. A full theory can be found in \cite{AN}. This multi-type model is useful in various fields, including biology, physics, and epidemiology, to study populations with different types of individuals, and how these types evolve over time. Recall that in the classical Galton-Watson process, there is only one type of individual, whereas in the multi-type case, there are $k$ different types: $1,2,\hdots,k$. We shall assume that the process starts with only one individual of a specific type initially, i.e. ${\bf Z}_0 = {\bf e}_i$. Let ${\bf Z}_n = (Z_{1n}, Z_{2n}, \ldots Z_{kn} )$ be the $k$-dimensional vector representing the number of individuals of each type in the $n$th generation of the branching process. Let ${\bf z} = (z_1, \hdots,z_k)$. For each type $i = 1, 2, \ldots, k$, denote by $p_i(\bf z)$ the probability that an individual of type $i$ has the offspring vector ${\bf z}$. Note that $\sum_{\bf z} p_i(\bf z) = 1$. These probabilities define the {\it offspring distribution} for individuals of type $i$.\medskip

To analyze the probabilistic structure of the offspring distribution, we can use generating functions. The probability generating function for a multi-type branching process with the offspring distribution $(p_i)$ is defined as
\be
\label{eqn:ProbGenFcn}
f_i( s_1,\hdots, s_k) = \sum_{{\bf z \in \mathbb{N}_0^k}} p_i({\bf z}) s_1^{z_1} s_2^{z_2} \cdots s_k^{z_k}, \,\,\,\,\,\,\,\,\,\, |s_1|, \hdots, |s_k| \leq 1.
\ee
We say $f_i$ is the {\it probability generating function} for the offspring distribution of individuals of type $i$.

For a $k$-type branching process, define the {\it matrix of means} as the $k \times k$ matrix ${\bf M} = [m_{ij}]$ such that
\be
\label{eqn:kTypeMeansMatrix}
m_{i j} = \mathbb{E}[Z_{j 1} \, | \, {\bf Z}_0 = {\bf e}_i].
\ee
In other words, the entry $m_{i j}$ represents the average number of offspring of type $j$ produced by an individual of type $i$ in one generation. In the case where the entries $m_{ij}$ vary, we have a multi-type  {\it inhomogenous} branching process.

\medskip

Let us now consider the {\it independent multi-type Poisson branching process}, a special case of the multi-type branching process where the number of offspring of each type of a single individual follow independent Poisson distributions. Then the entries of the corresponding mean matrix $M$ are the parameters of the Poisson distributions. This process has offspring distributions of the form
\be
\label{eqn:PoissonDist}
p_i({\bf z}) = \prod_{j=1}^k e^{-m_{ij}} \frac{(m_{ij})^{z_j}}{z_j !}.
\ee
Applying \eqref{eqn:PoissonDist} to \eqref{eqn:ProbGenFcn} and the independence of the number of offspring of
 each type of a single individual yield the following form for probability generating function of the multi-type Poisson branching process:
\be
\label{eqn:PoissonPGF}
f_i({\bf s}) = e^{\langle {\bf e}_i \, | {\bf M} \, | {\bf s} - \boldsymbol{1} \rangle}.
\ee
\vskip2mm

\subsection{Extinction Probability} 

As gelation is to vector multiplicative coalescent processes, a fundamental question in the study of branching processes is whether the process becomes {\it extinct} as the number of generations $n$ tends to infinity. Recall that in the classical case, if the expected number of the offspring distribution is less than (or equal to) 1, then eventually the population dies out. But if the expected number is greater than 1, then the probability of extinction is the smallest non-negative solution to the equation $s = \varphi(s)$, where  $\varphi$ is the probability generating function of the corresponding offspring distribution.\medskip

 For a multi-type branching process with offspring distribution $({\bf Z}_n)$, let ${\boldsymbol \eta} = (\eta_1, \ldots, \eta_k)$ denote the extinction probability vector, where $\eta_i$ denotes the extinction probability of an individual of type $i$ when the process begins with ${\bf Z}_0 = {\bf e}_i$. In this case, the phase transition for ${\boldsymbol \eta}$ is related to the spectral radius of the mean matrix ${\bf M}$, as stated in the following theorem (see \cite[Theorem 4.1]{AN}). A matrix ${\bf A}$ is {\it positive regular} if some power ${\bf A}^n$ has all entries strictly positive; a branching process is {\it singular} if each individual has exactly one child, and {\it nonsingular} otherwise.

\begin{thm}
\label{thm:GenExtinct}
Let $f_i({\bf s})$ be the probability generating function and ${\bf M}$ be the matrix of means of a positive regular and nonsingular multi-type branching process. Let $\rho({\bf M})$ denote the spectral radius of ${\bf M}$.  Furthermore, let ${\boldsymbol \eta}$ be the vector of extinction probabilities.

{\em (a)} If $\rho({\bf M})  \leq 1$, then ${\boldsymbol \eta} = {\bf 1}$, the $k$-dimensional vector $(1,1,\ldots,1)$.

{\em (b)} If $\rho({\bf M})   > 1$, then ${\boldsymbol \eta}$ is the smallest (by component) solution less than ${\bf 1}$ of the fixed point equations
 \be
\label{eqn:GenFixedPointEqns}
 f_i({\bf s}) = s_i \,,\,\, \forall\, i = 1, 2, \ldots, k.
\ee
\end{thm}

In the case of the multi-type Poisson branching process, the fixed point equations \eqref{eqn:GenFixedPointEqns} have the form
\be
\label{eqn:PoissonFixedPointEqns}
e^{\langle {\bf e}_i \, | {\bf M} \, | {\bf s} - \boldsymbol{1} \rangle} = s_i
\ee
which we will refer to as the {\it Poisson fixed point equations}.

\bigskip

\section{Correspondence between Gelation and Extinction Probability with Applications}

The bridge between the phase transitions of these two stochastic models can be seen in the equivalence of the multidimensional Lambert-Euler inversion equations \eqref{eqn:yzSolsCoord} and the Poisson fixed point equations \eqref{eqn:PoissonFixedPointEqns} with ${\bf M} = V D[{\boldsymbol\alpha}] t$.  Specifically,
\be
\label{eqn:EqnsEquiv}
e^{\langle {\bf e}_i \, | V D[{\boldsymbol\alpha}] t \, | {\bf s} - \boldsymbol{1} \rangle} = s_i \hsp \Longleftrightarrow \hsp y_i  e^{-\langle {\bf e}_i  | V | {\bf y} \rangle} = \alpha_i t  e^{-\langle {\bf e}_i  | V | \boldsymbol{\alpha} t \rangle}
\ee
when $s_i = y_i/(\alpha_i t)$. The identification ${\bf M} = V D[{\boldsymbol\alpha}] t$, or equivalently $m_{ij} = \alpha_j v_{ij} t$,
implies that the offspring distribution of the multi-type Poisson branching process \eqref{eqn:PoissonDist} has the form
\be
\label{eqn:NewPoissonDist}
p_i({\bf z}) = \prod_{j=1}^k e^{-\alpha_j v_{ij}t } \frac{(\alpha_j v_{ij}t )^{z_j}}{z_j !},
\ee
where ${\bf z} = (z_1,z_2, \dots,z_k).$

%

%
%

Let us note that the correspondence between the coalescent process $C$ and the branching process $B$ is not just through equations. There is an asymptotic coupling of $C$, $B$ and an inhomogeneous random graph $G$.
The connection between $C$ and $G$ is described in~\cite{A} and~\cite{KO}, and the connection between $G$ and $B$ was rigorously established in~\cite{BJR}. Here, we just outline the connections and make some remarks.

A given cluster $X$, in the coalescent process $C$, at time $t$, containing $m$ particles, was formed from $m-1$ distinct mergers of sub-clusters before time $t$. For the multiplicative kernel only, such mergers can be thought of as mergers between {\it individual particles} in the sub-clusters, and so together they correspond to a spanning tree $T_X$ among the particles in $X$. Such a spanning tree is illustrated in the first panel in Figure 2. $T_X$ corresponds to a connected component in an inhomogeneous random graph $G$. The part sizes and (matrix of) edge probabilities in $G$ correspond to the mass of particle types and particle interaction matrix of $C$: see Table 1 for the details. When $t$ is a constant, almost all the connected components of $G$ will be trees.

\vskip2mm
\begin{figure}[h]
    \centering
    \includegraphics[width=0.75\textwidth]{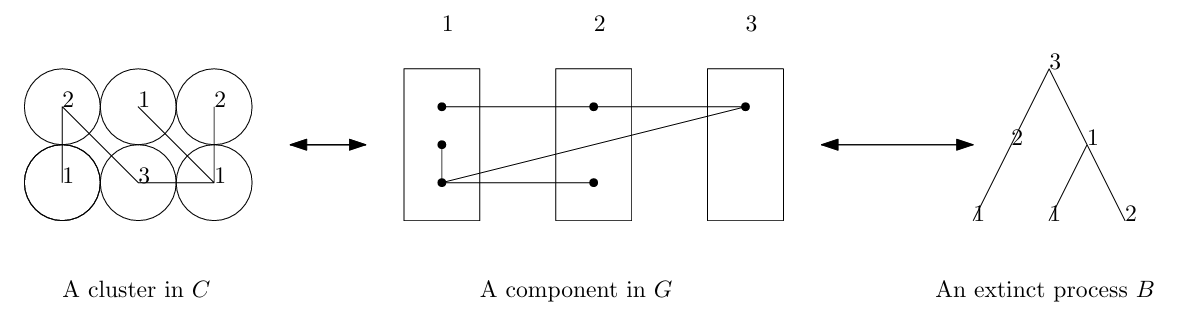}
    \caption{\small Three coupled models: $C$, $G$ and $B$}
    \label{fig:postge;}
\end{figure}

Next, we single out a vertex $v\in T_X$, of type $i$. (In Figure 2, $i=3$.) We can ``explore" the component of $G$ containing $v$ using a multi-type Poisson branching process $B$. In $B$, the parameters $t,\alpha_i$ and $v_{ij}$ have no direct interpretation, but they all factor into the means of the Poisson offspring distributions. Increasing $t$ increases the expected number $v_{ij}\alpha_jt$ of offspring of each type $j$, from each type $i$, in a linear fashion: so, as $t$ increases, the branching  process gets ``wider", not ``deeper". In particular, the branching process $B$ is {\it not} the time-reversed coalescent process $C$, the parameter $t$ is not the usual time parameter in the branching process, and the $n$th generation of $B$ does not correspond to the number of particles in $C$. However, the extinction probability $s_i(t)$ does correspond to the probability that a particle of type $i$ in $C$ has not joined the gel by time $t$.
\vskip2mm

For the branching process $B$, the time parameter has an interpretation in terms of the offspring distribution of a single parent. For a parent of type $i$, and offspring of type $j$, we construct a one-dimensional Poisson process $P$ on $\R^+$ of rate $\alpha_{ij}v_j$. For a given time $t$, the random number $N_{ij}(t)$ of offspring (of type $j$ produced by the parent of type $i$) is just the number of events in $P$ up to time $t$, i.e., in the interval $[0,t]$.


\begin{table}[h]\small
\[\begin{array}{|c|c|c|c|}\hline
{\rm Quantity}&{\rm Interpretation\ in\ }$C$&{\rm Interpretation\ in\ }$G$&{\rm Interpretation\ in\ }$B$\\\hline
t&{\rm time}&{\rm time}&{\rm see\ text}\\\hline
\alpha_i&{\rm proportion\ of\ particles\ of\ type\ }$i$&{\rm proportion\ of\ vertices\ of\ type\ }$i$&{\rm see\ text}\\\hline
v_{ij}&n\cdot{\rm bonding\ rate\ for\ types\ }i{\rm \ and\ }j&(n/t)\cdot{\rm \ edge\ probability\ for\ types\ }i{\rm \ and\ }j&{\rm see\ text}\\\hline
s_i(t)&\Prb({\rm particle\ of\ type\ }i{\rm \ not\ in\ gel})&\Prb({\rm vertex\ of\ type\ }i{\rm \ not\ in\ giant})&{\rm extinction\ probability\ for\ type\ }i\\\hline
v_{ij}\alpha_jt&-&-&\E({\rm type\ }j{\rm\ offspring\ from\ type\ }i)\\\hline
\alpha_is_i(t)&{\rm proportion\ of\ type\ }i{\rm \ not\ in\ gel}&{\rm proportion\ of\ type\ }i{\rm \ not\ in\ giant}&-\\\hline
\end{array}\]
\caption{The interpretation of various quantities in the different models}
\end{table}

Using this equivalence we get our first application that yields the gelation phenomenon for the vector multiplicative coalescent process discussed in Section \ref{GelationLE}.  We state the result again in the theorem below, followed by the short new proof that follows immediately from this equivalence.

\begin{thm}
Let $\zeta_{\bf x}(t)$ be the solution to the modified Smoluchowski equations, i.e., $\zeta_{\bf x}(t)$ represents the limiting fraction of clusters of size ${\bf x}$ in the vector multiplicative coalescent process. Then

{\em (a)} If $\rho(V D[{\boldsymbol\alpha}] t) \leq 1$, then the total mass $\sum_{\bf x}  \zeta_{\bf x}(t)|{\bf x}\rangle = |{\boldsymbol \alpha}\rangle $.

{\em (b)} If $\rho(V D[{\boldsymbol\alpha}] t)   > 1$, then the total mass $\sum_{\bf x}  \zeta_{\bf x}(t)|{\bf x}\rangle < |{\boldsymbol \alpha}\rangle $.

\vskip1mm
Consequently, the gelation time for the vector multiplicative coalescent process is $ T_{gel}={1 \over \rho ( V D[{\boldsymbol\alpha}] )}$.
\end{thm}
\bp
Applying the identification ${\bf M} = V D[{\boldsymbol\alpha}] t$ to Theorem \ref{thm:GenExtinct} yields

(a) If $\rho(V D[{\boldsymbol\alpha}] t)  \leq 1$, then ${\boldsymbol \eta} = {\bf 1}$, the $k$-dimensional vector $(1,1,\ldots,1)$.

{\em (b)} If $\rho(V D[{\boldsymbol\alpha}] t)   > 1$, then ${\boldsymbol \eta}$ is the smallest solution less than ${\bf 1}$ of the fixed point equations
\[ e^{\langle {\bf e}_i \, | V D[{\boldsymbol\alpha}] t \, | {\bf s} - \boldsymbol{1} \rangle} = s_i \]
Then the equivalence \eqref{eqn:EqnsEquiv} implies that $\boldsymbol{\alpha} t \boldsymbol{\eta} = {\bf y}$, which completes the proof.
\ep

\begin{remark}
In a recent work \cite{HKS}, Hoogendijk et al. described an alternative approach, which relies on establishing the equivalence between the Smoluchowski coagulation equations (see equation (\ref{eqn:SmolEqn})) and the inviscid Burgers equation using branching process techniques.
\end{remark}

The second application of the equivalence \eqref{eqn:EqnsEquiv} is a new series expression for the extinction probabilities of the Poisson branching process. From the equivalence \eqref{eqn:EqnsEquiv}, we have the extinction probability $\eta_\ell = y_\ell/(\alpha_\ell t)$, where $y_\ell$ is the smallest solution to the Lambert-Euler inversion equations.  Then by the total mass expression \eqref{eqn:TotalMass} and the explicit solution to the modified Smoluchowski equations \eqref{eqn:MSEsol}, we have
\be
\label{eqn:ExtProbSeries1}
\eta_\ell = \frac{1}{\alpha_\ell} \frac{y_\ell}{t} = \frac{1}{\alpha_\ell} \sum_{{\bf x}} \zeta_{{\bf x}}(t) x_\ell = \frac{1}{\alpha_\ell} \sum_{{\bf x}} {1 \over {\bf x}!}\boldsymbol{\alpha}^{\bf x} {\tau(K_k, x_i x_j v_{i,j}) \over {\bf x}^{\bf 1}}(V {\bf x})^{{\bf x}- {\bf 1}} e^{-\langle {\bf x} | V | \boldsymbol{\alpha} \rangle t} t^{\langle {\bf x} |\boldsymbol{1} \rangle -1} x_\ell.
\ee
From the relation ${\bf M} = V D[{\boldsymbol\alpha}] t$,
\[ \langle {\bf x} | V | \boldsymbol{\alpha} \rangle t = \langle {\bf x} | {\bf M} | \boldsymbol{1} \rangle \]
and
%

\[ (V {\bf x})^{{\bf x}- {\bf 1}} = \frac{1}{t^{\langle{\bf x | \bf 1 \rangle}- k}}({\bf M}D[\alpha_i^{-1}] {\bf x})^{{\bf x}-{\bf 1}}
\]
Moreover,
\[ \tau(K_k, x_i x_j v_{i,j}) = \frac{1}{t^{k-1}} \ \tau\left(K_k, x_i x_j \frac{m_{i,j}}{\alpha_j} \right). \]


Substituting the above three equations into \eqref{eqn:ExtProbSeries1} yields the following new result for the extinction probabilities of the multi-type Poisson branching process.
\begin{thm}
Suppose the matrix of means ${\bf M}$ for a multi-type Poisson branching process can be expressed as ${\bf M} = V D[{\boldsymbol\alpha}] t$ where $V$ is a nonnegative, irreducible and symmetric matrix, ${\boldsymbol \alpha} = (\alpha_1, \alpha_2, \ldots, \alpha_k) \in(0,\infty)^k$, and $t$ is a nonnegative real number.  Then the extinction probabilities of this process ${\boldsymbol \eta} = (\eta_1, \ldots, \eta_k)$ have the series form
\[ \eta_\ell = \frac{1}{\alpha_\ell} \sum_{{\bf x}} \frac{\boldsymbol{\alpha}^{\bf x}}{{\bf x}! \ {\bf x}^{\bf 1}} \ \tau\left(K_k, x_i x_j \frac{m_{i,j}}{\alpha_j} \right) ({\bf M}D[\alpha_i^{-1}] {\bf x})^{{\bf x}-{\bf 1}} e^{-\langle {\bf x} | {\bf M} | \boldsymbol{1} \rangle} x_\ell\,,\,\,\,\,\,\,\forall \ell = 1,2\ldots,k.\]
\end{thm}


In particular, for the homogenous multi-type Poisson branching process where $\boldsymbol{\alpha} = {\bf 1}$ and $V$ is the matrix with all entries $1$, the extinction probabilities have the form
\be
\label{eqn:HomExtProb}
\eta_\ell = \sum_{{\bf x}} \frac{1}{{\bf x}! \ {\bf x}^{\bf 1}} \ \tau\left(K_k, x_i x_j t \right) \langle {\bf x} | \boldsymbol{1} \rangle^{\langle {\bf x} | \boldsymbol{1} \rangle-k} t^{\langle {\bf x} | \boldsymbol{1} \rangle-k} e^{-\langle {\bf x} | \boldsymbol{1} \rangle kt} x_\ell\,,\,\,\,\,\,\,\forall \ell = 1,2\ldots,k.
\ee
Note that equation (\ref{eqn:HomExtProb}) is consistent with the single-type case. For example, consider the process with two types ${\bf x}=(x_1, x_2)$. Suppose $n= x_1 + x_2$. Then, (\ref{eqn:HomExtProb}) simplifies to
\begin{align*}
\eta_1  &= \sum_{x_1 \geq 1, n \geq 1} \frac{t}{x_1!\,(n-x_1))!} (tn)^{n - 2} e^{-2tn} x_1 \\
    & = \sum_{n \geq 1} \sum_{x_1=1}^{n-1} \frac{t}{x_1!\,(n-x_1)!} (tn)^{n - 2} e^{-2tn} x_1 \\
    & = \sum_{n \geq 1} t (tn)^{n - 2} e^{-2tn}  \left (\sum_{x_1=1}^{n-1} \frac{1}{x_1!\,(n-x_1))!} x_1 \right).\\
 \end{align*}
The inner sum can be rewritten as
$$\frac{1}{(n-1)!} \sum_{x_1=1}^{n-1} \frac{(n-1)!}{(x_1-1)!\,(n-x_1))!}$$
which by the binomial theorem, is equivalent to $\dfrac{2^{n-1}}{(n-1)!}$. Therefore,
\begin{align*}
    \eta_1 &= \sum_{n \geq 1} t (tn)^{n - 2} e^{-2tn} \dfrac{2^{n-1}}{(n-1)!}\\
    & =  \sum_{n \geq 1} \dfrac{n^{n - 1}}{n!} (2t)^{n-1}e^{-2tn}.
\end{align*}
Letting $z = 2te^{-2t}$, the last summation becomes
\begin{equation*}
    \eta_1 = \dfrac{1}{2t}  \sum_{n \geq 1} \dfrac{n^{n - 1}}{n!} z^n.
\end{equation*}
Let $T(z) = \sum_{n \geq 1} \dfrac{n^{n - 1}}{n!} z^n$, which is known as the {\it tree function} and it satisfies the functional equation $T(z) e^{-T(z)} = z$ and so, $T(z) = -W(-z)$, where $W$ is the {\it Lambert W function} (see \cite[section 2]{CGHJK}). Therefore, $ \eta_1  = -\frac{1}{2t} W(-2te^{-2t})$.\vskip2mm
On the other hand, for a single-type Poisson branching process with mean $\mu$, the extinction probablity $q$ is given by the fixed point equation
\begin{equation*}
     q= e^{\mu(q-1)}.
\end{equation*}
 Equivalently, $-\mu q = W(- \mu e^{-\mu})$, which implies $q = -\frac{1}{\mu}W(- \mu e^{-\mu})$. Hence, $\eta_1$ can be viewed as the extinction probability of a Poisson branching process with mean $\mu = 2t$, so if $\mu \leq 1$ (i.e. $t \leq 1/2$), then extinction is definite; $q=1$ (i.e. $\eta_1 = 1$).


\vspace{.5in}

\bibliographystyle{amsplain}

\begin{thebibliography}{10}

\bibitem{A} D. J.\ Aldous,
\newblock Deterministic and stochastic models for coalescence (aggregation and coagulation): a review of the mean-field theory for probabilists,
\newblock {\it Bernoulli} {\bf 6} (1999), 3–-48.

\bibitem{AS} N.\ Alon and J.\ Spencer,
\newblock The Probabilistic Method (3rd edition),
\newblock Wiley, New York (2008).

\bibitem{AIM} L.\ Andreis, T.\ Iyer and E.\ Magnanini,
\newblock Convergence of cluster coagulation dynamics,
\newblock arxiv preprint: 2406.12401, 2024.

\bibitem{AN} K. B.\ Athreya and P.\ E.\ Ney,
\newblock Branching Processes,
\newblock Springer-Verlag, Berlin (1972).

\bibitem{B}N.\ Berestycki,
\newblock Recent progress in coalescent theory,
\newblock {\it Ensaios Matem\'aticos} {\bf 16} (2009), 1--193.

\bibitem{Bol01} B. Bollob\'as,
\newblock Random Graphs (2nd edition),
\newblock Cambridge University Press (2001).

\bibitem{BJR} B.\ Bollob\'as, S.\ Janson and O.\ Riordan,
\newblock The phase transition in inhomogeneous random graphs,
\newblock {\it Random Structures \& Algorithms} {\bf 31} (2007), 3--122.

\bibitem{CGHJK}{R. M. Corless, G.H. Gonnet, D. E. Hare, D. J. Jeffrey and D. E. Knuth. On the Lambert W function, {\it Advances in Computational Mathematics}, 5(1) (1996), 329-359.}


\bibitem{ER} P.\ Erd\H{o}s and A.\ R\'enyi,
\newblock On the evolution of random graphs,
\newblock {\it Publ. Math. Inst. Hungar. Acad. Sci.} {\bf 5} (1960), 17--61.

\bibitem{EL} L. Euler,
\newblock De serie Lambertina Plurimisque eius insignibus proprietatibus,
\newblock {\it Acta Acad. Scient. Petropol.} {\bf 2} (1783), 29--51.

\bibitem{HP} D.\ Heydecker and R.I.\ Patterson,
\newblock Bilinear coagulation equations,
\newblock arxiv preprint 1902.07686 (2019).

\bibitem{HEZ} E. M.\ Hendriks, M.H.\ Ernst and R. M. Ziff,
\newblock Coagulation equations with gelation,
\newblock {\it Journal of Statistical Physics} {\bf 31} (1983), 519--563.

\bibitem{HKS} J.\ Hoogendijk, I.\ Kryven, and C.\ Schenone,
\newblock Gelation and localization in multicomponent coagulation with multiplicative kernel through branching processes,
\newblock {\it Journal of Statistical Physics} {\bf 191} (2024).

\bibitem{J}I.\ Jeon,
\newblock Existence of gelling solutions for coagulation-fragmentation equations,
\newblock {\it Comm. Math. Phys.} {\bf 194} (1998), 541--567.

\bibitem{K} J. F. C.\ Kingman,
\newblock The coalescent,
\newblock {\it Stochastic Processes and their Applications} {\bf 13(3)} (1982), 235--248.

\bibitem{KOY} Y.\ Kovchegov, P.\ T.\ Otto, and A.\ Yambartsev,
\newblock Cross-multiplicative coalescent processes and applications,
\newblock{\em  ALEA, Lat. Am. J. Probab. Math. Stat.} {\bf 18} (2021), 81--106.

\bibitem{KO} Y.\ Kovchegov and P.\ T.\ Otto,
\newblock Multidimensional Lambert-Euler inversion and vector-multiplicative coalescent processes,
\newblock{\em  Journal of Statistical Physics} {\bf 190} (2023).

\bibitem{JL} J.H. Lambert,
 \newblock Observationes variae in mathesin puram,
 \newblock {\it Acta Helveticae physico-mathematico-anatomico botanico-medica} {\bf 3} (1758), 128--168.

\bibitem{No} J.\ Norris,
\newblock Cluster coagulation,
\newblock {\it Comm Math Phys} {\bf 209} (2000), 407-–435.

\bibitem{R} F.\ Rezakhanlou,
\newblock  Gelation for Marcus-Lushnikov process,
\newblock  {\it Ann. Probab.} {\bf 41} (2013), 1806--1830.

\bibitem{S} M. V.\ Smoluchowski,
\newblock Drei vortrage uber diffusion, brownsche molekularbewegung und koagulation von kolloidteilchen,
\newblock {\it Z. Physik.} {\bf 17} (1916), 585--599.

\bibitem{VE}P.\ Van Dongen and M.\ Ernst,
\newblock On the occurrence of a gelation transition in Smoluchowski coagulation equation,
\newblock {\it J. Statist. Phys.} {\bf 44} (1986), 785--792.




\end{thebibliography}

\end{document}